\documentclass[a4paper,fleqn]{cas-dc}

\usepackage[numbers]{natbib}
\usepackage{lipsum} 
\usepackage{float}
\usepackage{subfig}

\def\tsc#1{\csdef{#1}{\textsc{\lowercase{#1}}\xspace}}
\tsc{WGM}
\tsc{QE}

\begin{document}
\let\WriteBookmarks\relax
\def\floatpagepagefraction{1}
\def\textpagefraction{.001}

\shorttitle{A generic decision support framework for retail inventory management}    

\shortauthors{H.J. Zietsman and J.H. van Vuuren}  

\title[mode = title]{A generic framework for decision support in retail inventory management}  



%

\author[]{H. Jurie Zietsman}[orcid=0000-0003-3683-0236]

\cormark[1]


\ead{20125844@sun.ac.za, uberziets@gmail.com, (+27) 76 278 9720}


\credit{Conceptualization, Methodology, Software, Validation, Data curation, Formal analysis, Writing --- original draft, Writing --- review \& editing}

\affiliation[]{organization={Stellenbosch Unit for Operations Research in Engineering, Department of Industrial Engineering, Stellenbosch University},
	addressline={Private Bag XI, Matieland}, 
	postcode={7602}, 
	country={South Africa}}

\author[]{Jan H. {van Vuuren}}[]


\ead{vuuren@sun.ac.za}

\ead[url]{http://www.vuuren.co.za/main.php}

\credit{Conceptualization, Methodology, Supervision, Resources, Writing --- review \& editing}





\begin{abstract}
The success of any business depends on how well it is able to satisfy customer demand, while remaining financially viable. Globalisation and the growth of e-commerce have resulted in retail businesses having to manage increasing numbers of products.  As the number of products sold increases, so does the complexity of the corresponding inventory management problem. High-quality decision making, in this regard, necessitates the utilisation of computerised decision support systems capable of accommodating the complexity presented by a large number of products. Existing frameworks for decision support in inventory management are mainly focused on a single facet of the inventory management problem --- generic, integrated frameworks seem to be absent from the literature. In this paper, a holistic framework design is proposed which integrates all of the major components expected to form part of a generic framework for inventory management. The objective of the framework is to provide decision support in respect of various inventory management operations, such as product segmentation, demand forecasting and determining the sizes and timings of replenishment orders, in pursuit of a desirable balance between the conflicting objectives in inventory management. 	
\end{abstract}

\begin{highlights}
\item A framework is proposed for decision support in inventory management
\item The framework is generic and may be adopted in any retail business
\item Results obtained during a verification analysis are competitive with benchmarks
\item The framework is demonstration within the context of a real-world case study
\item Results showed improvements over current operations of the case study business
\end{highlights}

\begin{keywords}
 Inventory management \sep Time series forecasting \sep Decision support systems \sep Supply chain management \sep
\end{keywords}

\maketitle

\section{Introduction}\label{Introduction}

Inventory management forms an integral part of supply chain management since it involves all the activities associated with balancing inventory levels and expected customer demand~\cite{Wild2017}. These activities may include planning, procurement, purchasing, warehousing, distribution and establishing policies. As such, inventory management plays a critical role in enabling a retail business to achieve an appropriate balance between satisfying demand while remaining financially sustainable. Achieving this balance necessitates supplying customers with sought-after products at reasonable prices within an acceptable time frame~\cite{Davis2016}.

The three most important pursuits of effective inventory management are to resolve the imbalance between supply and demand, minimise the risk of unavailability of \textit{stock keeping units}~(SKUs), and minimise costs associated with inventory management activities~\cite{Wild2017}. Rephrasing these motivators reveals three clear objectives, namely to (1) maximise customer service level, (2) minimise inventory on hand, and (3) minimise operating costs. The main decisions affecting any attempt at balancing between these objectives are (1) which SKUs have to be replenished, (2) when to place replenishment orders, and (3) the appropriate volumes of SKUs to request in these orders~\cite{Bushuev2015}. 

As the number of SKUs sold by a retailer increases, it becomes increasingly difficult to perform inventory management functions effectively. This is due to the complexity introduced by the aforementioned conflicting inventory management objectives. In order to make high-quality replenishment decisions in modern retail warehouses, it is therefore imperative that computerised \textit{decision support systems}~(DSSs) are employed to combat the complexity presented by large numbers of SKUs and varying customer demand.

Existing frameworks aimed at aiding in the decision making of inventory management practitioners are mainly focused on a single aspect of the inventory management problem, were developed for domain-specific problems, vaguely guide practitioners on the appropriate implementation of the framework, and lack flexibility in allowing user intervention and the utilisation of expert knowledge. 

In this paper, a generic framework for decision support in inventory management is presented which addresses the aforementioned shortcomings. The input to the framework is data related to the attributes, the historical demand, and the historical replenishment orders of SKUs, as well as information on suppliers. The framework facilitates the processing of these data, the clustering of SKUs based on their importance to the business, the forecasting of future demand for each SKU, and the recommendation of suitable inventory replenishment decisions per SKU (or group of SKUs), based on user-specified objectives. A DSS developed according to the proposed framework should aid practitioners in achieving a desired balance between satisfying customer demand and financial viability by providing good answers to the main decisions underlying sound inventory management, by synthesising recommended replenishment orders in a user-friendly manner.

The remainder of this paper is structured as follows: A brief overview of the relevant literature is provided in \S\ref{Literature_review}, focussing on existing inventory management frameworks. Our proposed framework is subsequently presented in some detail in \S\ref{IROS}. This is followed in \S\ref{Verification} by the verification of a computerised DSS instantiation of the framework, in respect of benchmark data from the literature, after which the practical applicability of the framework is showcased by means of a validation case study in \S\ref{Case_study}. The paper finally closes with a brief appraisal of its contributions and proposals for future avenues of research in \S\ref{Discussion} and \S\ref{Conclusions}, respectively.

\section{Literature review}\label{Literature_review}

The considerable body of research in the realm of inventory management may generally be partitioned according to the problem type considered. Prominent approaches in the literature towards three important inventory management operations are discussed briefly in \S\ref{lit_1}, after which the focus of the discourse shifts in \S\ref{lit_2} to existing frameworks for inventory management available in the literature. 

\subsection{Inventory management operations}\label{lit_1}

The practice of \textit{SKU classification} is aimed at tailored inventory management, production strategies, and demand forecasting per SKU group based on similar characteristics of SKUs. In a comprehensive review on SKU classification, four main attribute categories were identified which capture the majority of classification criteria present in literature~\cite{VanKampen2012}. These attribute categories are volume, product, customer, and timing. A multitude of approaches, both judgemental and statistical, have been employed to perform SKU classification. Judgemental techniques capture the expert knowledge of practitioners by, for example, estimating the \textit{criticality} of a product according to the \textit{vital, essential, desirable}~(VED) technique, or ranking SKU characteristics by means of the \textit{analytical hierarchy process}~(AHP) or the \textit{technique for order of preference by similarity to ideal solution}~(TOPSIS)~\cite{Bhattacharya2007}. The result forms part of the input required by statistical techniques from which final SKU clusters may be obtained. These statistical techniques include traditional ABC classification~\cite{Reid1987}, XYZ classification~\cite{Errasti2010}, graphical matrix approaches~\cite{DAlessandro2000}, neural networks\cite{Huiskonen2005}, optimisation models~\cite{HadiVencheh2010, Ng2007}, and even genetic algorithms\cite{Guvenir1998}.

The research field of \textit{time series forecasting} is well established and many extensive reference materials exist for both statistical forecasting models and machine learning forecasting models. Statistical models include naive models, moving averages, regression models, the class of exponential smoothing models, the class of \textit{autoregressive integrated moving average}~(ARIMA) models, and the class of \textit{autoregressive conditional heteroscedastic}~(ARCH) models, among others. Machine learning models, on the other hand, include tree-based models (such as the \textit{XGBoost} and \textit{LightGBM} algorithms), neural networks (such as multilayer perceptrons, recurrent neural networks~(RNNs), and convolutional neural networks~(CNNs)), and transformers, among others~\cite{Hyndman2020, Zhang2020}. 

Models specifically developed for forecasting \textit{intermittent time series} (which are commonly observed in the demand patterns of SKUs) have also been proposed. The first of these was the method of Croston~\cite{Croston1972}. A number of variations on Croston's method have subsequently been proposed. These include the \textit{Syntetos-Boylan approximation} (SBA)~\cite{Syntetos2005}, the \textit{Teunter-Syntetos-Babai} (TSB) model~\cite{Teunter2011}, and a modified SBA method~\cite{Babai2019}. Other models specifically developed for intermittent time series, namely the \textit{aggregate-disaggregate intermittent demand approach}\linebreak (ADIDA)~\cite{Nikolopoulos2011} and the \textit{multi-aggregation prediction algorithm}~(MAPA)~\cite{Kourentzes2014}, are based on temporal aggregation. Forecasting intermittent demand by means of neural networks has also been investigated~\cite{Kourentzes2013}.

\textit{Inventory replenishment models} are typically aimed at achieving some target service level, while minimising the related costs, which may include ordering costs, unit purchasing costs, holding costs and stockout costs~\cite{Grubbstroem1999}. Inventory replenishment models may broadly be classified as being either \textit{deterministic} or \textit{stochastic} in nature. In deterministic models, the underlying assumption is that the expected demand and lead times are known with certainty, while stochastic models account for some of the uncertainty present in real-world inventory management situations~\cite{Vandeput2020}. 

Within these two broad categories, inventory models may be classified along a number of dimensions, such as the assumed nature of customer demand (\textit{i.e.}\ static or dynamic), the number of SKUs considered (\textit{i.e.}\ single or multiple), and the number of echelons considered (\textit{i.e.}\ single or multiple)~\cite{Bushuev2015}. Another common differentiation between inventory model types is based on the nature of inventory monitoring, or the so-called \textit{review policy}~\cite{Davis2016}. Replenishment models are therefore also generally categorised as following either a \textit{continuous review} policy or a \textit{periodic review} policy. Finally, inventory models may also be classified based on the adopted planning horizon (\textit{i.e.}\ single-period or multi-period)~\cite{Jacobs2018}. 

The process of choosing specific models to employ in respect of different inventory management operations is not a trivial task, especially since an informed choice may result in better inventory control and increased profitability as a result of improved product availability, increased sales, and a reduction of inventory levels~\cite{Wild2017}. 

\subsection{Existing inventory management frameworks}\label{lit_2}

The existing literature on frameworks for inventory management also reveals a clear need for providing adequate decision support in multiple areas of inventory management. This need is highlighted by the large number of frameworks or DSSs available, as well as the range of inventory management areas for which they have been developed. 

These frameworks primarily differ in terms of (1) the assumed replenishment system characteristics (\textit{e.g.}\ the nature of demand, the number of SKUs accommodated, and the number of echelons considered), (2) the targeted inventory management function (\textit{e.g.}\ SKU classification, forecasting, the generation of replenishment policies, or a combination of these), and (3) the application domain (\textit{e.g.}\ focused on a particular business or industry, or generally applicable). In many of the framework descriptions an emphasis is placed on utilising the abundance of data available to businesses. 

\citeauthor{Kobbacy1999}~\cite{Kobbacy1999} were among the first to propose a framework for an intelligent inventory management system in 1999. They attempted to bridge the gap between the theory on inventory models and the practical application thereof, highlighted as early as 1981 by \citeauthor{Silver1981}~\cite{Silver1981}. A number of the existing frameworks are specifically aimed at determining appropriate inventory replenishment policies based on a particular method of SKU segmentation~\cite{Errasti2010, Kobbacy1999,Spyridakos2009, Wanke2014}. In these cases, frameworks have typically been designed around a single or multiple specific SKU segmentation metrics, the most common of which is the nature of customer demand for an SKU. The focus of the framework proposed in 2014 by \citeauthor{Wanke2014} \cite{Wanke2014}, for example, was on low-consumption SKUs. 

The framework proposed in 2016 by Krishnadevarajan \textit{et al.}\ \cite{Krishnadevarajan2016} presented a structured approach towards selecting metrics to employ when classifying SKUs. A multi-criteria approach towards SKU classification is important, since the historical customer demand pattern associated with an SKU is not necessarily the only factor on which inventory replenishment policies in a business are or should be based.

The purpose of the framework proposed by El-Aal \textit{et al.}\ \cite{ElAal2010} in 2010 was to evaluate and compare different inventory control policies based on any number of evaluation metrics. In 2015, Van den Berg \textit{et al.}\ \cite{Berg2013} provided a framework for determining which key performance indicators to measure at different stages of the inventory management process. The combined application of these frameworks may prove beneficial to a business in terms of calculating combined key performance indicators for inventory management.

A few domain-specific frameworks or systems have also been proposed~\cite{Shang2008, Spyridakos2009}. Future applications of instances of these frameworks are, however, limited due to their domain-specific nature. Furthermore, the framework proposed in 2009 by Spyridakos \textit{et al.}\ \cite{Spyridakos2009} is the only framework reviewed which takes maintenance into account (due to its domain-specific development).

The frameworks proposed in 2020 by Errasti \textit{et al.}~\cite{Errasti2010} and in 2011 by \citeauthor{Cadavid2011}~\cite{Cadavid2011} are respectively more holistic in nature and attempt to integrate SKU classification, demand forecasting, and inventory replenishment models. The domain-specific DSS proposed by Shang \textit{et al.}\ \cite{Shang2008} in 2008 and the framework of Aviv~\cite{Aviv2003}, dating from 2003, are the only other frameworks which accommodate the forecasting of demand. With the increasing availability of large volumes of sales and inventory-related data, the ability to forecast demand may prove to be an indispensable component of any modern inventory management framework.

A holistic framework, in this sense, refers to a framework that integrates, in some manner, all of the aforementioned components expected to be present in a generic inventory management framework. No such a generic, integrated framework is available in the literature, to the best of our knowledge. Moreover, the existing frameworks do not allow for the combination of SKUs in replenishment orders in pursuit of economies of scale. The design of the framework proposed in this paper is distinguished from these frameworks by its \textit{integration} of the following features: 
\begin{enumerate}
	\item A multi-criteria clustering approach towards SKU segmentation in pursuit of prioritising certain SKUs in terms of resource allocation.
	\item A tailored time series forecasting capability aimed at SKU demand estimation which allows for the possibility of time series clustering, as well as the ensembling forecasting models.
	\item A robust order replenishment model capable of effectively handling exceptions in customer demand and consolidating SKUs into inventory replenishment batches pursuant of economies of scale. 
\end{enumerate}

\section{The proposed framework}\label{IROS}

In this section, we present in some detail the \textit{inventory replenishment operations support}~(IROS) framework. An overview of its design is first provided, and this is followed by descriptions of the primary framework components.

\subsection{The framework architecture}\label{Framework_overview}

The IROS framework is designed with the main objective of \textit{unifying} the existing inventory replenishment approaches in the literature, resulting in a framework which is \textit{generalisable} across industries and business types, as well as sufficiently flexible to adapt to the specific inventory replenishment needs of a particular business. Moreover, the framework is aimed at being \textit{accessible} to a wide range of users, since its implementation facilitates models ranging from simple to state-of-the-art. The framework is also designed to be \textit{interactive}, in the sense that user input is elicited from processes within the framework.

\begin{figure*}[h] 
	\centering		
	\includegraphics[scale=0.75]{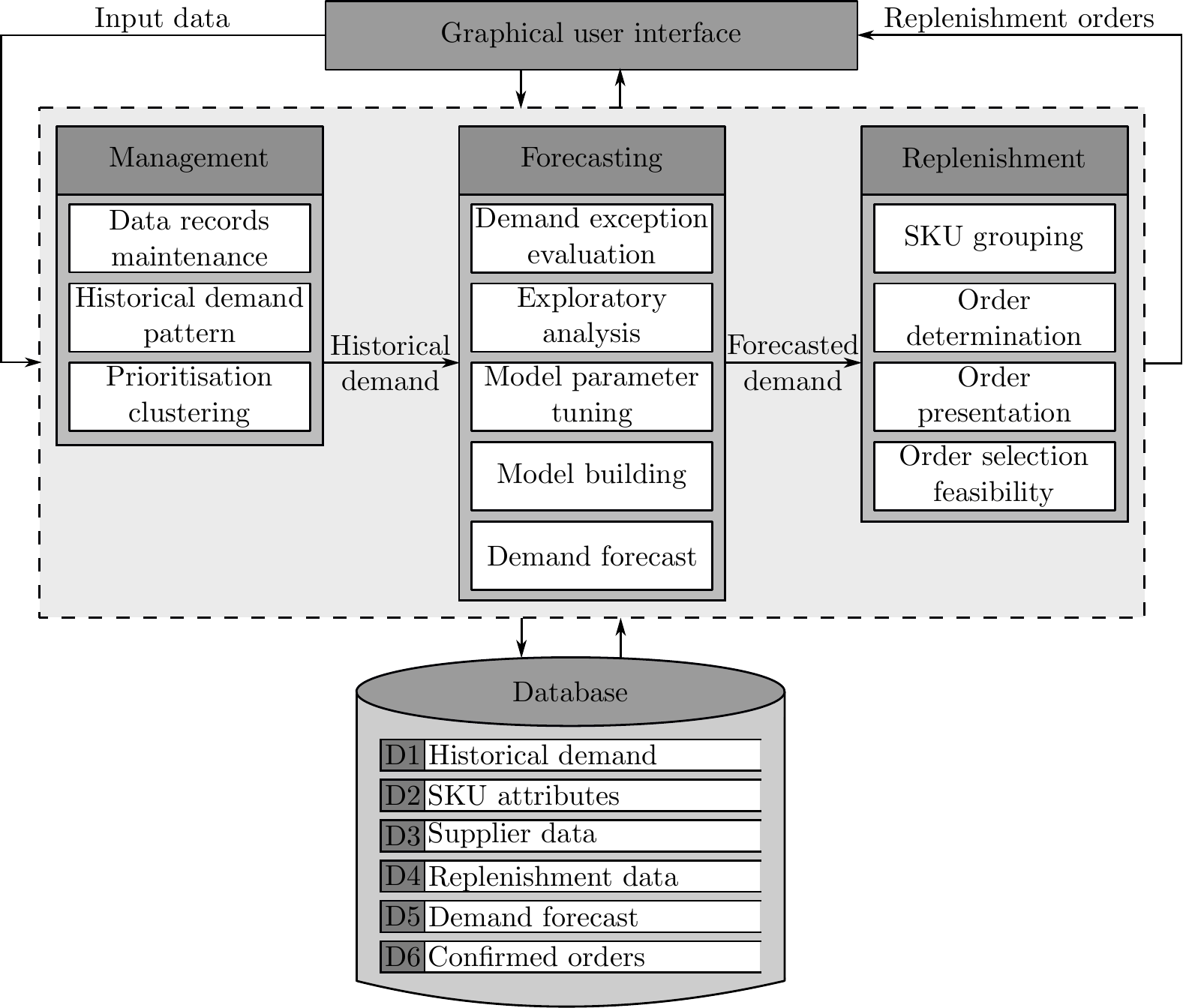}
	\caption[A high-level overview of the newly proposed IROS framework]{A schematic high-level overview of the proposed IROS framework architecture, indicating the underlying modules of the primary components.}
	\label{fig:Framework_Overview}
\end{figure*}

A high-level schematic overview of the framework architecture is presented in Figure~\ref{fig:Framework_Overview}. The framework comprises three primary components, and the functions which form part of each are described with reference to the component's constituent modules (or processes). These modules form the data manipulation and computational backbone of the IROS framework and may collectively be referred to as the \textit{central processing} component. The modules are generic in nature, and allow users to implement algorithms or models freely which satisfy their particular subjective preferences. 

The user provides input data \textit{via} a graphical user interface and receives decision support in the form of recommended replenishment orders as output. Additional inputs are also elicited from the user during the execution of each of the primary components. Such interactions are necessary due to the generic nature of the framework design. For instance, the user has to specify performance measures which are most applicable to the particular business. The database forms the foundation of the IROS framework, as it stores the data required by all of the framework components and facilitates communication between the component modules. 

The framework is designed to follow a \textit{transient} strategy in which the planning horizon for inventory replenishment continually shifts one time step into the future in the form of a moving time window. It is assumed that historical sales of SKUs are sufficiently (although not necessarily solely) indicative of customer demand and may be employed to forecast future customer demand. Moreover, a \textit{demand-driven} inventory replenishment strategy is followed --- that is, the replenishment component modules are aimed at satisfying forecast customer demand. It is also assumed that SKU attribute data are available which detail at least the current inventory level, current number of back-orders, \textit{minimum order quantity}~(MOQ), unit mass, unit dimensions, and the supplier of each SKU. Data on the typical lead times of SKUs or suppliers are also assumed to be available.

The primary components of the framework and their respective constituent modules are discussed in greater detail in the remainder of this section. The discussion is aided by the use of \textit{data flow diagrams}~(DFDs).

\subsection{The management component}\label{Management_component}

\begin{figure*}[h] 
	\centering
	\includegraphics[scale=1.5]{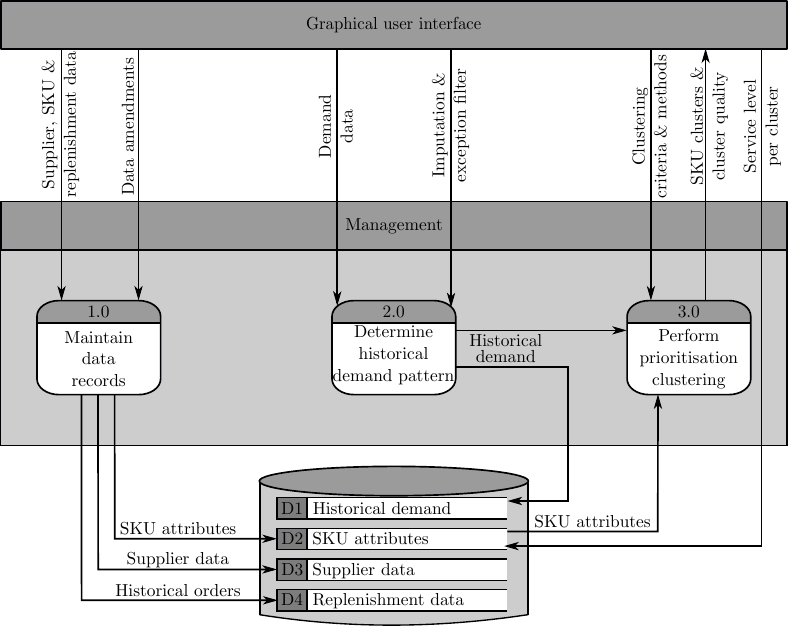}
	\caption[Level-one DFD of the framework management component]{Level-one DFD of the management component.}
	\label{fig:Framework_Management}
\end{figure*}

The management component of the IROS framework comprises three modules numbered 1.0 to 3.0, as shown in its level-one DFD in Figure~\ref{fig:Framework_Management}. This component may be viewed as a preprocessing step to the subsequent components. The purpose of Module~1.0 is to manage the data records and entails inputting and updating records aimed at ensuring that the most accurate data are available. Module~2.0 is designed to determine a historical demand pattern for each SKU, based on past customer demand, and includes the imputation of missing values, as well as the filtering of exceptions in customer demand. In Module~3.0, SKUs are clustered into groups based on their importance to the business. 

The input provided by the user to Module~1.0 consists of multiple data sets, which describe the SKU attributes, suppliers, and historical replenishment orders of the business. The information captured in the features of these data sets may differ, depending on the requirements of the chosen models and approaches employed in the remainder of the framework. SKU attributes would typically contain information on the unit cost, unit sales price, unit mass, unit dimensions, supplier(s), current inventory level, and current number of back-ordered units, for each unique SKU. 

The lead time associated with ordering an SKU may also be included under SKU attributes, for example, or as part of supplier data where it would indicate the lead time associated with all of the SKUs ordered from a particular supplier. The supplier data set may also contain information stipulating the weight capacity and/or volume capacity of containers or trucks used to transport orders from each supplier. The features of these data sets are not fixed, but some data are expected to remain relatively constant over time, such as the SKU attributes and supplier data. Replenishment data provide information on the historical order patterns of the business for each SKU. 

These data records have to be updated when, for instance, new replenishment orders are placed. Amendments to current data records are also required when, for example, the MOQ of an SKU is adjusted by a supplier. The SKU attributes and supplier data stores also have to be updated in the case where new SKUs are added to the product offering. 

In Module~2.0, the historical data related to customer demand have to be consolidated so as to represent the true demand for each SKU as accurately as possible. In cases where stockouts occurred and customer demand could not be fulfilled, or when returns are made due to an oversupply of SKUs, historical sales may, however, not be entirely reflective of the true customer demand. Quotes which are not accepted by customers may also be included when determining historical customer demand, if the assumption holds that the business is applying the correct pricing to its products.

Additional data preparation may have to be performed on the consolidated demand data, such as imputing missing values in cases where there are known data omissions. In cases where there was simply no customer demand for an SKU, however, data imputation would be erroneous. Data preparation may also include transforming the consolidated data by filtering out exceptions in the historical demand that the user does not wish to account for during inventory management. Note the use of the term \textit{exception} instead of \textit{outlier}, since the aim of this filtering process is exactly that --- to remove demand exceptions with known underlying reasons which are not expected to recur. The models employed as part of the forecasting component should be capable of dealing with outliers appropriately, as well as with time series exhibiting intermittent demand patterns. 

A prioritisation clustering procedure is performed in Module~3.0, based on one or multiple clustering criteria (\textit{i.e.}\ SKU attributes) and a specified clustering method. The objective of this segmentation is to facilitate inventory management improvement by creating clusters of SKUs for which similar service levels may be specified by the user. By adopting such a prescriptive service level perspective towards replenishment operations, the IROS framework follows the so-called inventory optimisation approach~\cite{Davis2016}. This approach is also echoed in the design of the forecasting and replenishment components.

The user should select a subset of clustering criteria which best describe the priorities of the business in terms of resource allocation. The clustering method employed may either be judgemental or statistical in nature, or a sequential combination of both may be applied. The user receives a summary of the clusters obtained and a measure of cluster quality as output, and may alter the clustering criteria or method if desired. The service level specified per SKU cluster is finally stored for later use in the replenishment component.

After being invoked initially in the IROS framework, the management component modules are not necessarily required during each subsequent iteration of its execution, but rather only when required by the user. 

\subsection{The forecasting component}\label{Forecasting_component}

\begin{figure*}[h] 
	\centering
	\includegraphics[scale=1.5]{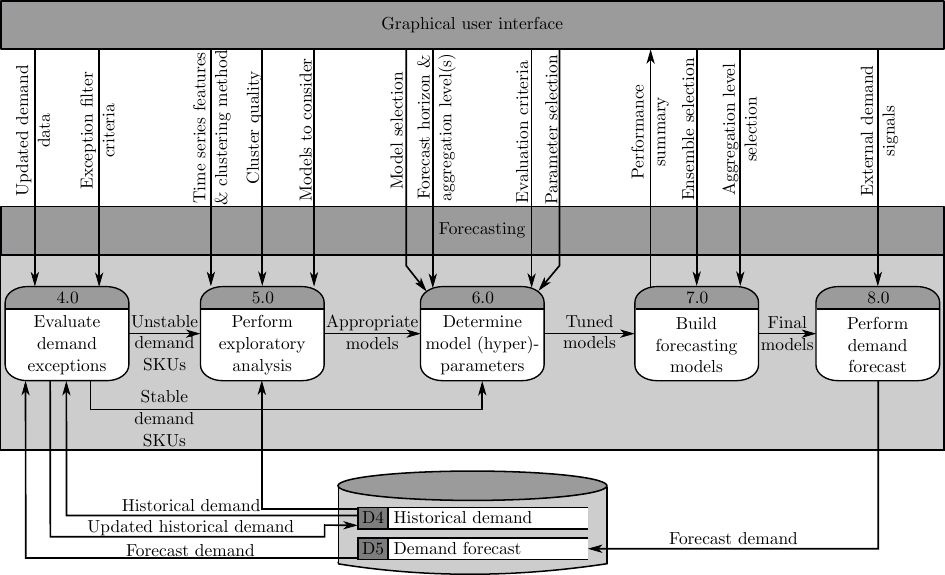}
	\caption[Level-one DFD of the framework forecasting component]{Level-one DFD of the forecasting component.}
	\label{fig:Framework_Forecast1}
\end{figure*}

The second component of the IROS framework, its forecasting component, comprises five modules numbered 4.0 to 8.0, as shown in its level-one DFD in Figure~\ref{fig:Framework_Forecast1}. The primary aim in this component is to obtain an accurate customer demand forecast for each SKU based on the historical demand pattern determined in Module~2.0 of the management component.

Module~4.0 is designed to evaluate the difference between forecast demand and actual customer demand, and to report on any significant demand exceptions. After having filtered out demand exceptions (as in Module~2.0), the historical demand for each SKU is updated accordingly. An exploratory analysis is then performed in Module~5.0 on the historical demand patterns of SKUs exhibiting unstable demand patterns. The benefit of such an exploratory analysis lies in determining appropriate forecasting models for individual time series, based on clusters of similar time series. 

The (hyper)parameters of appropriate forecasting models are subsequently tuned in Module~6.0 for each SKU, based on its historical demand and certain inputs provided by the user. The model building process in Module~7.0 allows the user to combine the tuned models from Module~6.0 by a process of model ensembling, delivering the final forecasting models to be employed per SKU. The expected demand for each SKU is finally predicted in module~8.0.

Focussing on the unique case of the initial iteration of the IROS framework execution, Module~4.0 is bypassed and an exploratory analysis is performed on the historical demand of each SKU. This consists of extracting relevant time series features which describe the characteristics of the demand time series. Based on the characteristics of each time series, certain forecasting models may be more appropriate than others. Different time series may therefore be clustered together, based on the similarity of their characteristics. Appropriate forecasting models for each cluster of time series may then be determined. The purpose of this procedure is to lessen the computational burden of the remaining modules in the forecasting component, as only the appropriate models (per cluster) are included in the following computations. Module~5.0 may alternatively be disregarded by allowing each SKU to form part of its own cluster. 

In all other cases (after the initial framework iteration), the historical demand for each SKU have to be updated in Module~4.0 to reflect the most recent customer demand data, in line with the transient replenishment strategy adopted in the framework. This is achieved by evaluating the difference between the demand forecast and the actual customer demand. Any significant differences in the form of anomalies or sustained deviations are reported to the user who may, analogously to the exception filtering performed in Module~2.0, opt to include or exclude demand exceptions during inventory replenishment planning, based on expert knowledge. During this process, SKUs exhibiting significant customer demand exceptions (not filtered out by the user) are labelled as \textit{unstable}. The process of exploratory analysis is repeated for such unstable demand SKUs, as the time series characteristics of their historical demand may have changed. More importantly, the appropriate forecasting models considered in Modules~6.0--8.0 may also have changed. SKUs exhibiting no significant demand exceptions are labelled as \textit{stable} and do not undergo another round of exploratory analysis.

The (hyper)parameters of appropriate forecasting models have to be tuned separately for each SKU in Module~6.0. Appropriate forecasting models are therefore fitted on historical data and the best parameters identified for each model per SKU are employed during the subsequent model building process in Module~7.0. The user may select additional models to include or possibly exclude models deemed inappropriate, based on expert knowledge.  

Important factors which may influence the quality of the final demand forecast are the forecast horizon over which models should be applied (\textit{e.g.}\ one year or six months), the frequency of these forecasts (\textit{e.g.}\ daily, weekly or monthly), and the corresponding level of aggregation of historical demand data (\textit{e.g.}\ daily, weekly or monthly). The manner in which data should be partitioned into a forecast training set and validation set are also required inputs by the user. Model performance is subsequently evaluated based on a user-specified evaluation criterion. The \textit{mean absolute scaled error}~(MASE) is arguably the most appropriate, since it accommodates all forecasting situations, all forecasting methods, and all time series types~\cite{Hyndman2006}.

During the model building process in Module~7.0, the tuned models are combined according to some method of ensembling, such as a simple mean or a weighted aggregation of model predictions. The user then receives a performance summary of the best-performing forecasting models per SKU, according to an appropriate evaluation criterion.  Based on this performance summary, the user may confirm the final forecasting model selection for each SKU, as well as the selected aggregation level. The final model for each SKU may comprise an ensemble of models or an individual model, depending on their performance and user preferences. A cross-validation approach towards model evaluation in both Modules~6.0 and 7.0 is recommended in order to yield models which generalise well.

The purpose of Module~8.0 is to generate the final customer demand forecast to be used in the replenishment component of the IROS framework. During this process, additional information which may influence customer demand may be specified by the user. Such customer demand signals represent external events, such as planned SKU promotions, and may take the form of an anticipated percentage increase in the demand for a particular SKU over some period of time. After having accounted for any external demand signals, the final demand forecast is performed and stored for each SKU.

\subsection{The replenishment component}\label{Replenishment_component}

\begin{figure*}[h] 
	\centering
	\includegraphics[scale=1.5]{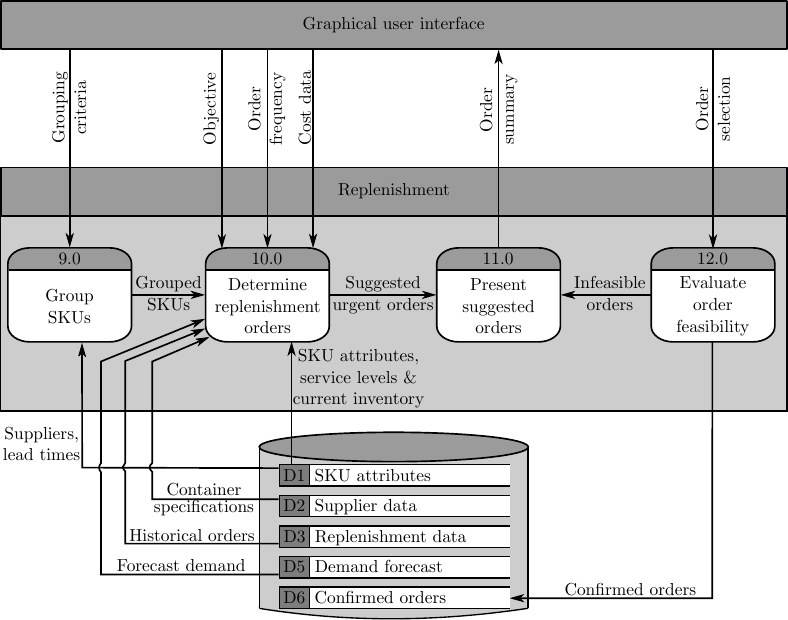}
	\caption[Level-one DFD of the framework replenishment component]{Level-one DFD of the replenishment component.}
	\label{fig:Framework_Replenishment1}
\end{figure*}

The working of the final component of the IROS framework, its replenishment component, is illustrated in its level-one DFD in Figure~\ref{fig:Framework_Replenishment1}. It comprises Modules~9.0 to 12.0, which are designed to provide the user with recommended urgent replenishment orders that satisfy the target service level (established in the management component) in respect of the forecast customer demand (obtained from the forecasting component). 

Module~9.0 is responsible for creating groups of SKUs based on the possibility of these SKUs being consolidated into the same replenishment orders. These grouped SKUs undergo a temporal batching procedure in Module~10.0, aimed at determining a good inventory replenishment strategy per group of SKUs, for the entire planning horizon.  The suggested urgent replenishment orders returned as output by Module~10.0 are passed to Module~11.0, where a summary of recommended urgent orders is generated and presented to the user. The user may then alter the details of the suggested orders, based on expert knowledge, upon which the resulting order selection is evaluated for feasibility in Module~12.0 in terms of the MOQ of each SKU  and restrictions on the volume and weight of delivery containers, as imposed by suppliers. Based on this evaluation, orders are either confirmed or again presented to the user as being infeasible. In such a case, another order selection is elicited.

The motivation behind the grouping performed in Module~9.0 lies in the flexibility it imparts on the temporal batching procedure of Module~10.0. Grouping SKUs based on whether they may be consolidated together in the same order (\textit{e.g.}\ SKUs ordered from the same supplier or location) allows for the problem of determining replenishment orders to be solved per group of SKUs. Such an approach may achieve a reduction in inventory costs in a manner analogous to achieving \textit{economies of scale}. The user may specify the grouping criteria employed in Module~9.0. 

Thereafter, the process in Module~10.0 is executed for each group of SKUs. The purpose of this module is to generate a combined order replenishment strategy for the group and, consequently, an individualised replenishment strategy for each SKU in the group. This goal is also in line with the principles of inventory optimisation~\cite{Davis2016}. The process should take into account the forecast demand, SKU attributes (such as unit cost, unit sales price, unit mass, unit dimensions, and MOQ), current inventory level and number of back-ordered units of each SKU, as well as historical orders made by the business, and volumetric and weight restrictions of orders imposed by suppliers. The forecast customer demand, as well as the lead times associated with replenishment orders may include an uncertainty measure in order to better account for stochasticity. 

The user is also required to supply relevant cost data if such data have not been stored in the database. The desired order frequency (\textit{i.e.}\ how often orders are placed) is also a required input. The chosen order frequency should be considered in conjunction with the final aggregation level specified in Module~7.0 of the forecasting component.

The resulting replenishment strategy should attempt to satisfy the desired customer service level of each SKU, as mentioned before, but is furthermore evaluated based on a user-specified objective considered to be desirable for the business, such as minimising cost or maximising profit. A good replenishment strategy should therefore achieve an acceptable balance between the conflicting objectives of inventory management discussed in Section~\ref{Introduction}.

The generic nature of the IROS framework allows for the adoption of a variety of approaches towards the implementation of Module~10.0. The adopted approach should be based on the user-specified inventory management objective. Standard inventory replenishment models may be employed. The adopted approach determines the complexity and required computational capacity of the process. The output returned by Module~10.0 aids the decision making process by providing (near-)optimal recommended replenishment orders in an attempt to answer the main questions in respect of inventory management, by specifying when replenishment orders should be placed, as well as which SKUs and the quantity of each to include in replenishment orders.

\begin{figure}[pos=h] 
	\centering
	\includegraphics[scale=0.8]{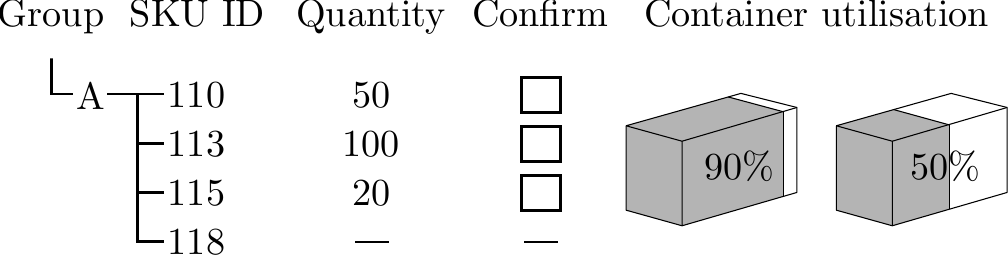}
	\caption{An example of presenting recommended replenishment orders.}
	\label{fig:Suggested_Orders}
\end{figure}

The suggested urgent orders are passed to Module~11.0 where they are aggregated and displayed in a user-friendly and intuitive order summary, an example of which is shown in Figure~\ref{fig:Suggested_Orders} for a single SKU group. The user may then alter the suggested orders based on expert knowledge. The final order selection of the user is evaluated for feasibility in Module~12.0, if any changes are made (\textit{e.g.}\ do the specified number of units still fit within the volume and/or weight restrictions of the supplier's delivery container or truck). Orders that remain feasible after user selection are labelled as \textit{confirmed} orders. In the case where the order selection by the user is not feasible, the relevant orders are again displayed to the user \textit{via} Module~11.0. This iterative process is repeated until all of the current urgent orders have been confirmed.

Execution of the forecasting and replenishment components are repeated according to the time steps included by the order frequency of the replenishment component. This transient approach therefore entails updating the historical demand to reflect the most recent customer demand and generating a new demand forecast. The demand-chasing order replenishment strategy is subsequently redetermined based on the new demand forecast and updated inputs. The IROS framework is ultimately aimed at supporting the development of an inventory replenishment DSS tailored to the specific requirements of a particular business. 

\section{Verification analyses}\label{Verification}

In an attempt to showcase the functionality of the framework components, a DSS instantiation developed according to the proposed framework was implemented in Python~3.7 on a personal computer and verified against benchmark data available in the literature. The modular design of the framework allows for the verification of the components of the DSS instantiation to be conducted on an individual basis. For each benchmark data set considered, the relative performances of the models included in the DSS instantiation were compared with benchmark performances documented in previous studies.

\subsection{Exception detection}\label{Exceptions}

The publicly available \textit{Yahoo Webscope S5} data set was employed when verifying the implementation of the exception detection process. The data comprises a total of 367 time series, for which the true anomalous observations are known. These series are partitioned into four smaller data sets called A1, A2, A3 and A4. Subset A1 contains real-world data of aggregated user logins to the Yahoo network, while the time series in the other subsets were produced synthetically. A variety of methods embedded in the \texttt{anomaly detection toolkit}~(ADTK) Python library~\cite{ADTK} were incorporated into the DSS instantiation to detect exceptions in the data. The performances of the DSS instantiation methods, as well as the reported performances of detection methods in the literature (taken as benchmarks), are provided in Table~\ref{tab:Exceptions} based on the mean \textit{area under curve} (AUC) the \textit{receiver operation characteristics} (ROC) curve~\cite{Skiena2017} for each data set. 

\begin{table}[pos=h] 
	\caption{Mean AUC scores achieved by anomaly detection methods, with the best performances typeset in boldface for each data subset.}
	\label{tab:Exceptions}
	\begin{tabular}{lccccc}
		& \multicolumn{4}{c}{Data subset} &  \\
		\multicolumn{1}{c}{Model} & A1 & A2 & A3 & A4 & Mean \\ \hline
		IQR & 1.000 & 0.500 & 0.500 & 0.500 & 0.625 \\
		Quantile & 0.982 & 0.976 & 0.726 & 0.180 & 0.823 \\
		AR & \textbf{0.976} & 0.749 & 0.997 & 0.998 & 0.93 \\
		Seasonal & 0.862 & 0.998 & \textbf{0.999} & \textbf{1.000} & 0.965 \\ \hline
		Twitter~\cite{Twitter2015} & 0.824 & 0.500 & 0.618 & 0.653 & 0.649 \\
		iForest~\cite{Liu2008} & 0.888 & 0.662 & 0.628 & 0.633 & 0.703 \\
		DeepAnt~\cite{Munir2018} & 0.898 & 0.961 & 0.928 & 0.86 & 0.912 \\
		FuseAD~\cite{Munir2019} & 0.947 & \textbf{0.999} & \textbf{0.999} & 0.966 & 0.978 \\ \hline
	\end{tabular}%
\end{table}

The Twitter benchmark model was developed with the objective of automatic anomaly detection in a practical and robust manner~\cite{Twitter2015}. The \textit{isolation forest}~(iForest) model is an ensemble of isolation trees in which anomalies are `isolated' from the non-anomalous observations of the time series~\cite{Liu2008}. Two state-of-the-art unsupervised models are also included as benchmarks. In the DeepAnt model, proposed by Munir \textit{et al.}\ \cite{Munir2018} in 2018, a CNN is fitted to the time series before anomalous residuals are detected. The FuseAD model, proposed by Munir \textit{et al.}\ \cite{Munir2019} in 2019, combines statistical and deep learning models by first fitting an ARIMA model to the original time series, after which a CNN corrects the residuals, before anomalies are detected.

A comparison of the performances of the detection methods implemented in the DSS instantiation with those of the benchmark methods revealed increasing AUC scores with an increase in the complexity of the detection models. The two best-performing DSS detection methods, the \textit{autoregressive}~(AR) and seasonal detectors, were able to outperform all but one of the benchmark models, on average.  

Considering the relative simplicity of the DSS methods, especially when compared with the CNN-based approaches of the DeepAnt and FuseAD models, they may be judged to have performed better than expected. The simpler methods require much lower volumes of data to train and, in the context of inventory management, their implementation may prove advantageous since the volume of available data may be limited and inventory managers are anticipated to be capable of interpreting and adjusting such methods more easily. The results achieved by the exception detection methods of the DSS instantiation were deemed satisfactory in respect of the benchmark models.

\subsection{SKU classification}

Clustering methods from the \texttt{scikit-learn}~\cite{scikit-learn} and \texttt{scikit-learn som}~\cite{sklearn_som} Python libraries were adopted for both SKU classification and time series clustering. A benchmark data set for multi-criteria inventory classification, introduced by Reid~\cite{Reid1987}, is employed in an attempt to compare the relative performances of the different SKU segmentation methods embedded in the DSS instantiation with published benchmark performances by a number of authors~\cite{HadiVencheh2010, Ng2007}. The data set consists of 47 SKUs with three features each, namely \textit{annual sales value}, \textit{average unit cost}, and \textit{lead time}, listed here in order of importance. Each of the three SKU attributes is also assumed to be positively correlated to the importance of an SKU from the view point of a retail business. The variable values for each SKU are scaled to the range $[0,1]$ in accordance with the convention of previous authors. 

In previous studies, the SKUs in the data set were partitioned into three clusters and each SKU was allocated to one of a number of clusters based on pre-determined occupation proportions, such that Cluster A contains 10 SKUs, Cluster B contains 14 SKUs, and Cluster C contains 23 SKUs. In the DSS instantiation, the most appropriate number of clusters according to the cluster inertia for the $k$-means algorithm and the dendrogram produced for agglomerative hierarchical clustering was also determined to be three.

The quality of the clusters obtained by these two methods, as well as by \textit{self-organising map} (SOM) clustering and two benchmark models, were evaluated in respect of the silhouette coefficient~\cite{Kaufman1990}, the Davies Bouldin (DB) index~\cite{Davies1979}, and the Cali\'{n}ski-Harabasz (CH) index~\cite{Calinski1974}. The results are provided in Table~\ref{tab:cluster_separation}. The Ng and Hadi-Vencheh optimisation models achieved a higher cluster quality according to the silhouette score and the DB index, while the DSS clustering methods performed the best in respect of the CH index. Based on these metrics, there is no single method that obtained significantly better clusters than any other method.

\begin{table}[pos=h]
	\caption{Cluster separation for the different models employed during SKU prioritisation clustering, with the best model performances typeset in boldface.}
	\label{tab:cluster_separation}
	\begin{tabular}{lccc}
		& \multicolumn{3}{c}{Cluster quality metric} \\
		Model & Silhouette & DB index & CH index \\ \hline
		$k$-means & 0.11 & 4.086 & 88.618 \\
		Hierarchical & 0.054 & 2.656 & \textbf{89.160} \\
		SOM & 0.054 & 2.656 & \textbf{89.160} \\ \hline
		Ng~\cite{Ng2007} & \textbf{0.114} & 1.756 & 32.674 \\
		Hadi-Vencheh~\cite{HadiVencheh2010} & \textbf{0.114} & \textbf{1.714} & 22.69 \\ \hline
	\end{tabular}%
\end{table}

In the absence of cluster proportions being provided as inputs, the traditional clustering methods assigned nine SKUs to Cluster~A, 22 SKUs to Cluster~B, and 16 SKUs to cluster~C (based on the majority vote) --- with the majority of differences occurring in the intersection between Clusters~B and~C. Based on these results, it was concluded that the prioritisation clustering process implemented in the DSS instantiation achieved satisfactory results in respect of the benchmark data set when compared with published results.

\subsection{Time series clustering}

The extraction of time series features in the DSS instantiation was implemented by utilising the \texttt{darts}~\cite{darts}, \texttt{nolds}~\cite{nolds}, \texttt{scipy}~\cite{scipy}, and \texttt{statsmodels}~\cite{statsmodels} Python libraries. The features included are indicated in Table~\ref{tab:ts_features}, as well as whether they were calculated for the original time series and/or a trend-and-seasonally adjusted version of the time series.

\begin{table}[pos=h] 
	\caption{The features extracted for time series clustering~\cite{Wang2006}.}
	\label{tab:ts_features}
	\centering
	\begin{tabular}{lcc}
		\multicolumn{1}{c}{Feature} & \begin{tabular}[c]{@{}c@{}}Original \\ time series\end{tabular} & \begin{tabular}[c]{@{}c@{}}Adjusted \\ time series\end{tabular} \\ \hline
		Trend & --- & \checkmark \\
		Seasonality & --- & \checkmark \\
		Serial correlation & \checkmark & \checkmark \\
		Skewness & \checkmark & \checkmark \\
		Kurtosis & \checkmark & \checkmark \\
		Stationarity & \checkmark & \checkmark \\
		Chaos & \checkmark & --- \\
		Self-similarity & \checkmark & --- \\
		Variation & \checkmark & --- \\
		Intermittency & \checkmark & --- \\ \hline
	\end{tabular}
\end{table}

The implementation of the time series clustering process was verified in respect of a data set from the \textit{University of California, Riverside}~(UCR) time series classification archive~\cite{UCRArchive2018}. The \textit{MixedShapesRegularTrain} was used, comprising 500 time series. This data set was selected for verification of the time series clustering process since the time series are labelled as belonging to one of five classes, depending on the characteristics of the time series. 

\textit{Principal component analysis}~(PCA) was applied to the feature space and, based on the resulting variance drop-off, the clustering procedure was restricted to the top three principal components. Although it is known that there are five distinct classes of time series in the data set, the most appropriate number of clusters for the DSS clustering methods was investigated. The largest cluster separation was indeed achieved when five clusters were employed. The clusters obtained by the $k$-means algorithm achieved the best performance in respect of the silhouette score, DB index, and CH index and were taken as the final time series clusters. The performance of the time series clustering procedure was subsequently evaluated in respect of the labels of the time series originating from the original benchmark data set. 

The clusters obtained from the $k$-means algorithm were labelled according to the five class labels in the original data. The number of series belonging to each of the benchmark classes, for each cluster, is shown in Table~\ref{tab:mode_cluster}. Inspecting the most prevalent $k$-means cluster for each time series class revealed the following mappings: Cluster~4 $\mapsto$ Classes~1 and~2, Cluster~2 $\mapsto$ Classes~3 and~4, and Cluster~3 $\mapsto$ Class~5, while Clusters~1 and~5 were not the most prevalent members for any of the time series classes.

\begin{table}[pos=h]
	\caption{The number of cluster members per time series class.}
	\label{tab:mode_cluster}
	\begin{tabular}{ccccccc}
		&  & \multicolumn{5}{c}{$k$-means cluster} \\
		&  & 1 & 2 & 3 & 4 & 5 \\ \cline{3-7} 
		\multirow{5}{*}{\rotatebox[origin=c]{90}{Class}} & \multicolumn{1}{c|}{1} & 17 & 10 & 23 & \textbf{32} & \multicolumn{1}{c|}{18} \\
		& \multicolumn{1}{c|}{2} & 0 & 2 & 28 & \textbf{69} & \multicolumn{1}{c|}{1} \\
		& \multicolumn{1}{c|}{3} & 5 & \textbf{78} & 7 & 0 & \multicolumn{1}{c|}{10} \\
		& \multicolumn{1}{c|}{4} & 3 & \textbf{64} & 7 & 3 & \multicolumn{1}{c|}{23} \\
		& \multicolumn{1}{c|}{5} & 0 & 12 & \textbf{86} & 0 & \multicolumn{1}{c|}{2} \\ \cline{3-7} 
	\end{tabular}
\end{table}

In order to quantify the performance of the time series clustering process, a confusion matrix was constructed. It is assumed that the above cluster to class mappings are correct and that Classes~1 and~2 are combined, as are Classes~3 and~4. The resulting confusion matrix is given in Table~\ref{tab:confusion_cluster}. 

\begin{table}[pos=h]
	\caption[Confusion matrix for time series clustering verification]{Confusion matrix for the verification of time series clustering.}
	\label{tab:confusion_cluster}
	\begin{tabular}{ccccc}
		&  & \multicolumn{3}{c}{Predicted class} \\
		&  & 1 \& 2 & 3 \& 4 & 5 \\ \cline{3-5} 
		\multirow{3}{*}{\rotatebox[origin=c]{90}{True class}} & \multicolumn{1}{c|}{1 \& 2} & \textbf{101} & 12 & \multicolumn{1}{c|}{51} \\
		& \multicolumn{1}{c|}{3 \& 4} & 3 & \textbf{142} & \multicolumn{1}{c|}{14} \\
		& \multicolumn{1}{c|}{5} & 0 & 12 & \multicolumn{1}{c|}{\textbf{86}} \\ \cline{3-5} 
	\end{tabular}
\end{table}

From this confusion matrix, it is clear that 50.5\% of the time series in Classes~1 and~2 were correctly clustered together, 71\% of Classes~3 and~4, and 86\% of Class 5. The time series in Classes~1 and~2 were the most difficult to cluster, leading to the largest number of incorrect classifications. It was, however, concluded that the performance of the time series clustering procedure is satisfactory, based on the manner in which distinguishing features between time series exhibiting different characteristics could be identified and clustered upon. This conclusion seems plausible, especially in view of demand forecasting being the ultimate aim. Time series representing demand patterns typically exhibit more prominent differences in characteristics (such as the difference between intermittent series and continuous series). 

\subsection{Forecasting}\label{Verification_forecasting}

The verification of the DSS forecasting processes was conducted on the well-known M4 forecasting competition data~\cite{Makridakis2018b}. The M4 data consists of $100\,000$ time series labelled according to frequency (\textit{i.e.}\ yearly (Y), quarterly (Q), monthly (M), weekly (W), daily (D) and hourly (H)) and domain~\cite{Makridakis2020}. An attempt was made to determine whether the performance of a particular forecasting model on a subset of time series from a time series cluster could be taken as a good proxy for the performance of the model on all the time series in the cluster. The motivation was the potential for reducing computational requirements when evaluating relative model performances for the time series in a cluster. 

A stratified subset of the M4 data was selected in order to evaluate the validity of such an approach. A total of 360 time series were randomly selected from each of the sets corresponding to the six frequencies (except the weekly frequency due to the availability of only 359 time series). These were also selected to be representative of the six domains. The composition of the data subset is given in Table~\ref{tab:M4_Data2}.

\begin{table}[pos=h] 
	\centering
	\caption{The subset of the M4 data employed during verification.}
	\label{tab:M4_Data2}
	\begin{tabular}{lrrrrrrr}
		Domain & Y & Q & M & W & D & H & Total \\ \hline
		Micro & 60 & 60 & 60 & 112 & 70 & 0 & 362 \\
		Industry & 60 & 60 & 60 & 6 & 70 & 0 & 256 \\
		Macro & 60 & 60 & 60 & 41 & 70 & 0 & 291 \\
		Finance & 60 & 60 & 60 & 164 & 70 & 0 & 414 \\
		Demo & 60 & 60 & 60 & 24 & 10 & 0 & 214 \\
		Other & 60 & 60 & 60 & 12 & 70 & 360 & 622 \\ \hline
		Total & 360 & 360 & 360 & 359 & 360 & 360 & 2\,159 \\ \hline
	\end{tabular}%
\end{table}

Time series clustering was applied to this subset based on the features listed in Table~\ref{tab:ts_features}. After performing PCA, the top five principal components were adopted, and three time series clusters was deemed the most appropriate, based on cluster separation. The $k$-means algorithm again performed the best. An analysis revealed that the majority of time series belonging to the same frequency of the M4 data were also grouped together in the same time series cluster. The three clusters could be mapped to the different frequencies according to the mappings shown in Table~\ref{tab:M4_mappings}. 

\begin{table}[pos=h]
	\caption{Cluster mappings and best-performing forecasting models for the M4 data subset.}
	\label{tab:M4_mappings}
	\begin{tabular}{cllll}
		Cluster & \multicolumn{1}{c}{Frequency} & \multicolumn{3}{c}{Best-performing models} \\ \hline
		1 & W, D  & Naive S & ARIMA & ETS \\
		2 & Y, Q, M & Theta & Theta\_bc & Naive S \\
		3 & H & LR & RF & LightGBM \\ \hline
	\end{tabular}
\end{table}

The forecasting models  implemented in the DSS instantiation included the benchmark models proposed by the organisers of the M4 competition~\cite{Makridakis2018b}, as well as a CNN, tree-based models (a random forest~(RF), XGBoost, and LightGBM), and a linear regression~(LR) model. The machine learning models were implemented based on the time series reduction approach for applying machine learning models to time series data. The \texttt{sktime} and \texttt{sktime\_dl} Python libraries were employed for the implementation of these models. The relative performances of these models were subsequently evaluated by means of 5-fold rolling origin cross validation performed on the training data (based on the \textit{overall weighted average}~(OWA) metric of the M4 competition). The three best-performing models per cluster are given in Table~\ref{tab:M4_mappings} in descending order.

Models were ranked based on their \textit{median} cross-\linebreak validated OWA scores for each of the clusters. One of the best-performing models per cluster, and associatively per frequency, coincided with the best-performing model for each frequency across the \textit{full} M4 data set (out of those considered), as documented in the literature~\cite{Loening2020}. Based on this finding, it was concluded that a subset of a time series cluster (based on the features listed in Table~\ref{tab:ts_features}) was indeed sufficiently representative of all the time series in that cluster in the context of determining the most probable best-performing forecasting model(s) for the time series of the cluster.

The best-performing models per time series cluster were then passed as inputs to the model building phase (\textit{i.e.}\ model ensembling). During the model building phase, all possible model ensembles of the three models were evaluated in respect of the \textit{full} M4 data set, again based on 5-fold rolling origin cross validation. The best-performing individual model or model ensemble for each time series was then employed to perform the final forecast, which was evaluated against the test set. A comparison between the performance of the forecasting process of the DSS instantiation, aggregated for all frequencies, and the top three models submitted to the M4 competition, is provided in Table~\ref{tab:M4_Performance}. Based on the results obtained, the DSS instantiation would have achieved a ranking of 7\textsuperscript{th} out of the top 50 original competition entries, based on OWA.

\begin{table}[pos=h] 
	\caption{Benchmark comparison on M4 competition data~\cite{Makridakis2020}.}
	\label{tab:M4_Performance}
	\begin{tabular}{lllcc}
		Author/model & sMAPE & MASE & OWA & Rank \\ \hline
		Smyl & 11.374 & 1.536 & 0.821 & 1 \\
		Montero-Manso, \textit{et al.} & 11.720 & 1.551 & 0.838 & 2 \\
		Pawlikowski, \textit{et al.} & 11.845 & 1.547 & 0.841 & 3 \\
		IROS & 11.928 & 1.595 & 0.856 & 7 \\ \hline
	\end{tabular}%
\end{table}

\begin{figure}[pos=h] 
	\centering
	\includegraphics[scale=0.47]{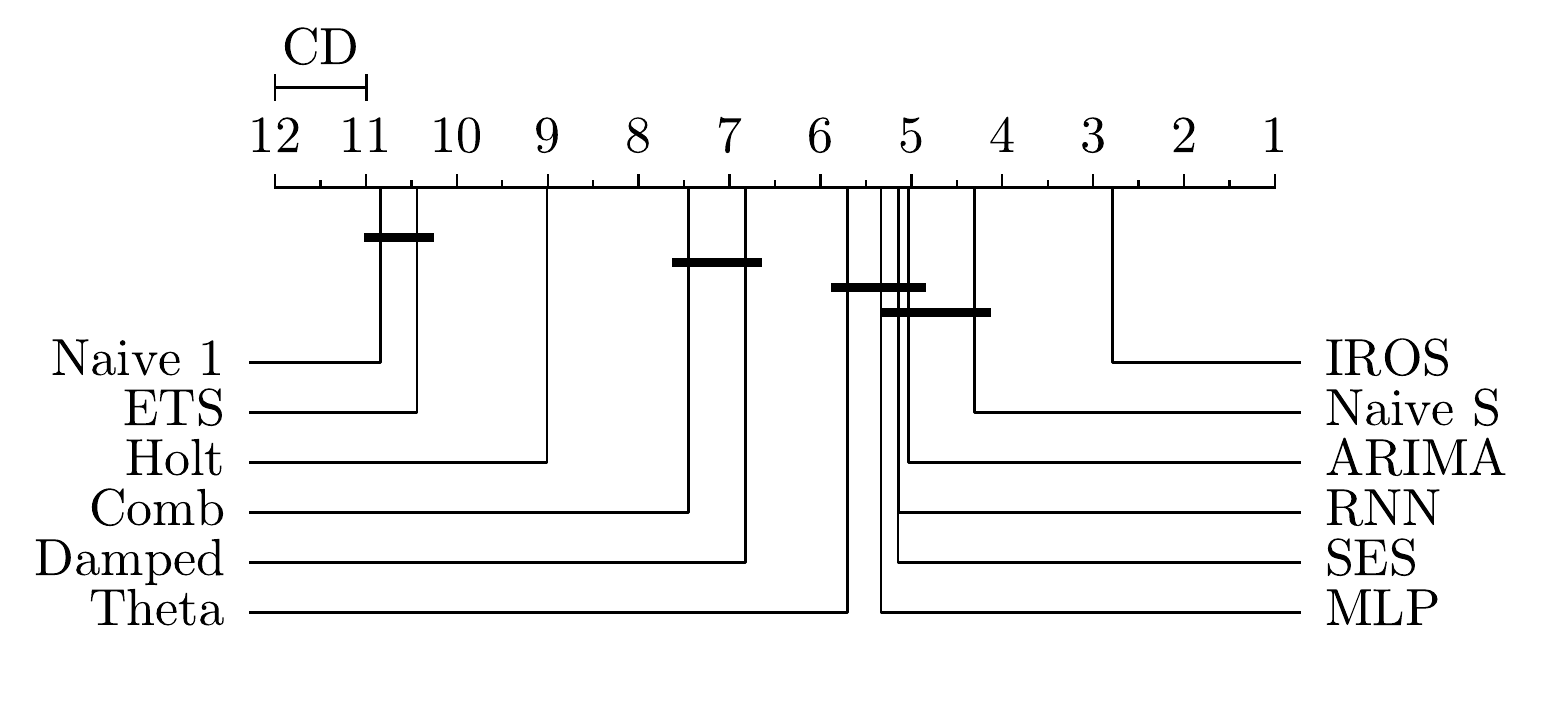}
	\caption{CD diagram comparing the DSS instantiation to the M4 benchmark models for the hourly time series data.}
	\label{fig:M4_Hourly_CD}
\end{figure}

It was furthermore determined whether the performance differences between the DSS instantiation and the benchmark models proposed for the M4 competition were statistically significant. For this evaluation, the OWA was calculated individually for each time series. The null hypothesis of the Friedman test was rejected for all frequencies at a level of $\alpha=0.01$ statistical significance. The Nemenyi post-hoc test was utilised to infer which performance differences were statistically significant. The DSS instantiation performed well in respect of all six frequencies, ranking first for the Y, M and H frequencies, second for the Q and W frequencies, and fourth for the D frequency.  Only for the H frequency was the difference in performance between the DSS instantiation and \textit{all} of the benchmark models statistically significant. The \textit{critical distance}~(CD) diagram for this frequency may be found in Figure~\ref{fig:M4_Hourly_CD}. 

Based on the models employed during the model building phase, these results are consistent with those reported by other authors~\cite{Loening2020}. Based on the performance of the DSS instantiation in comparison with individual benchmark models and the top-performing models submitted to the M4 forecasting competition, it was concluded the DSS instantiation achieved satisfactory results.

\section{Validation case study}\label{Case_study}

In an attempt to demonstrate the practical applicability of the IROS framework, a validation case study was undertaken in which the verified DSS instantiation was applied to real-world data of a retail business operating in South Africa. The business is a leading distributor of fasteners to the shipping, automotive, construction, mining, transport, agriculture, and manufacturing industries, among others. The sales of the business vary between approximately 400\,000 and 600\,000 transactions annually. Five financial years of data were made available for the case study. The first four years formed the \textit{training} set, while the final year formed the \textit{comparison} set.

\subsection{Historical demand patterns}\label{Historical_patterns}

\begin{figure}[pos=h]
	\centering
	\subfloat[Weekly aggregation]{
		\includegraphics[width = 3.6cm]{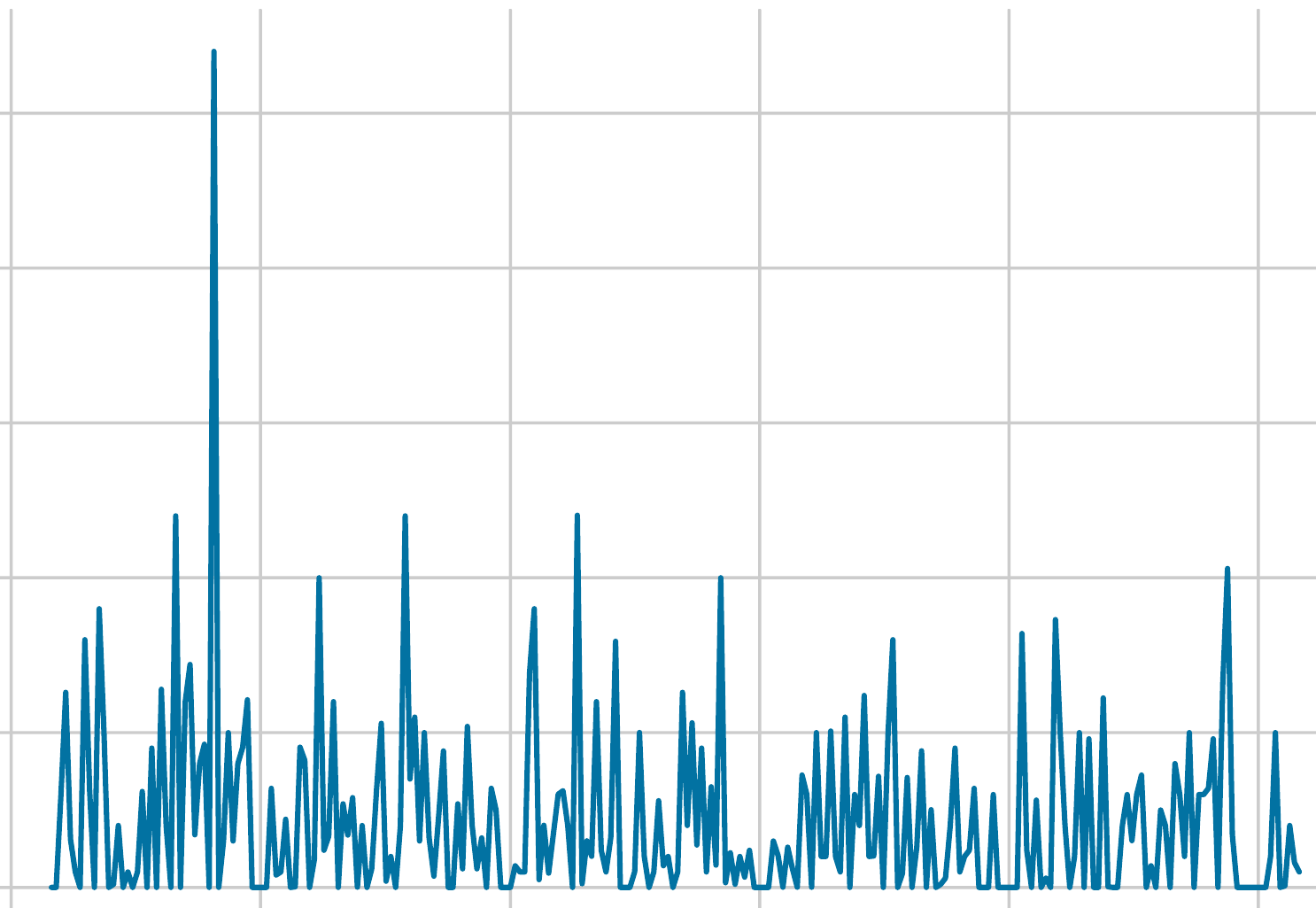}
		\label{fig:ts2}
	}
	\quad
	\subfloat[Monthly aggregation]{
		\includegraphics[width = 3.6cm]{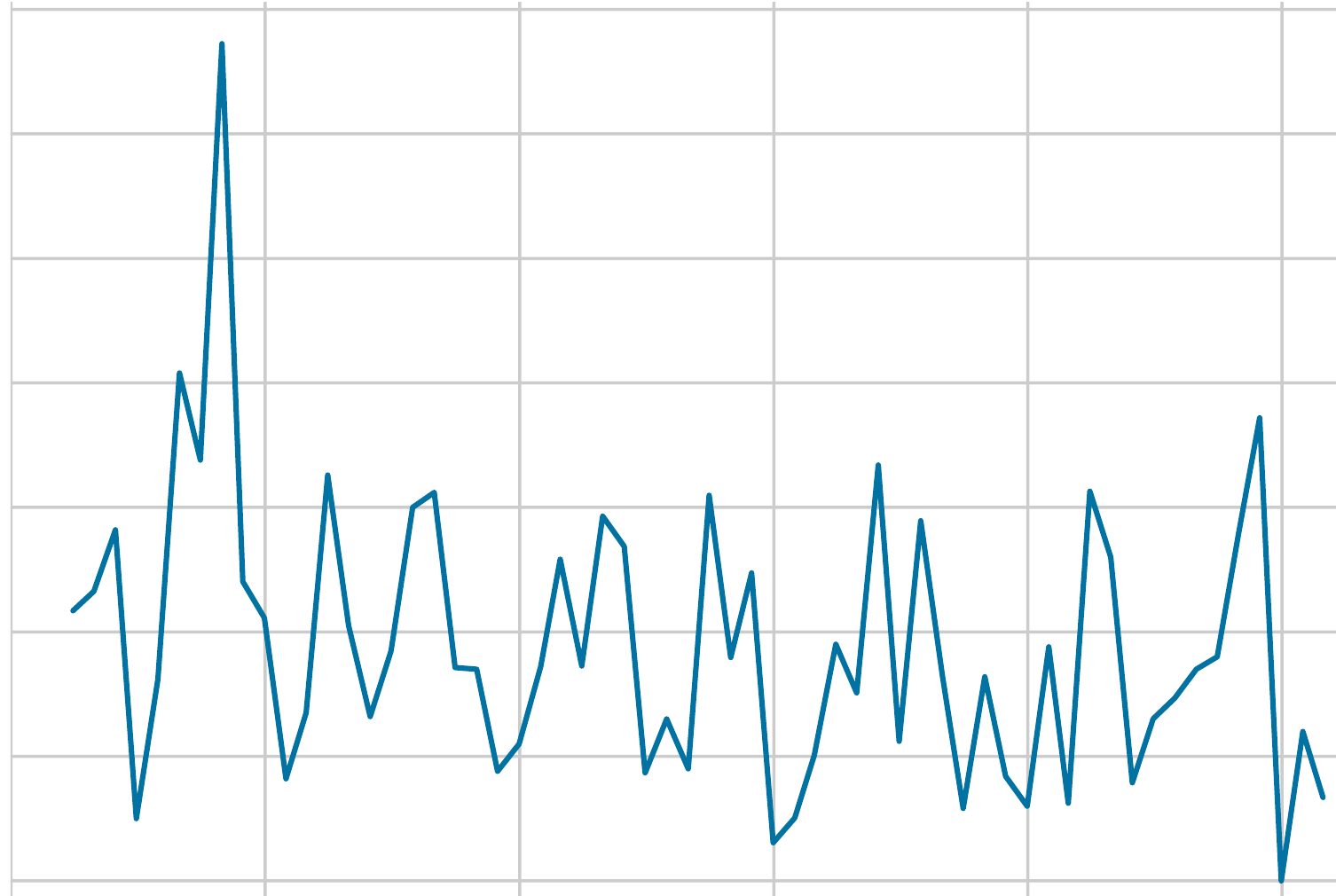}
		\label{fig:ts3}
	}
	\quad
	\subfloat[Weekly aggregation]{
		\includegraphics[width = 3.6cm]{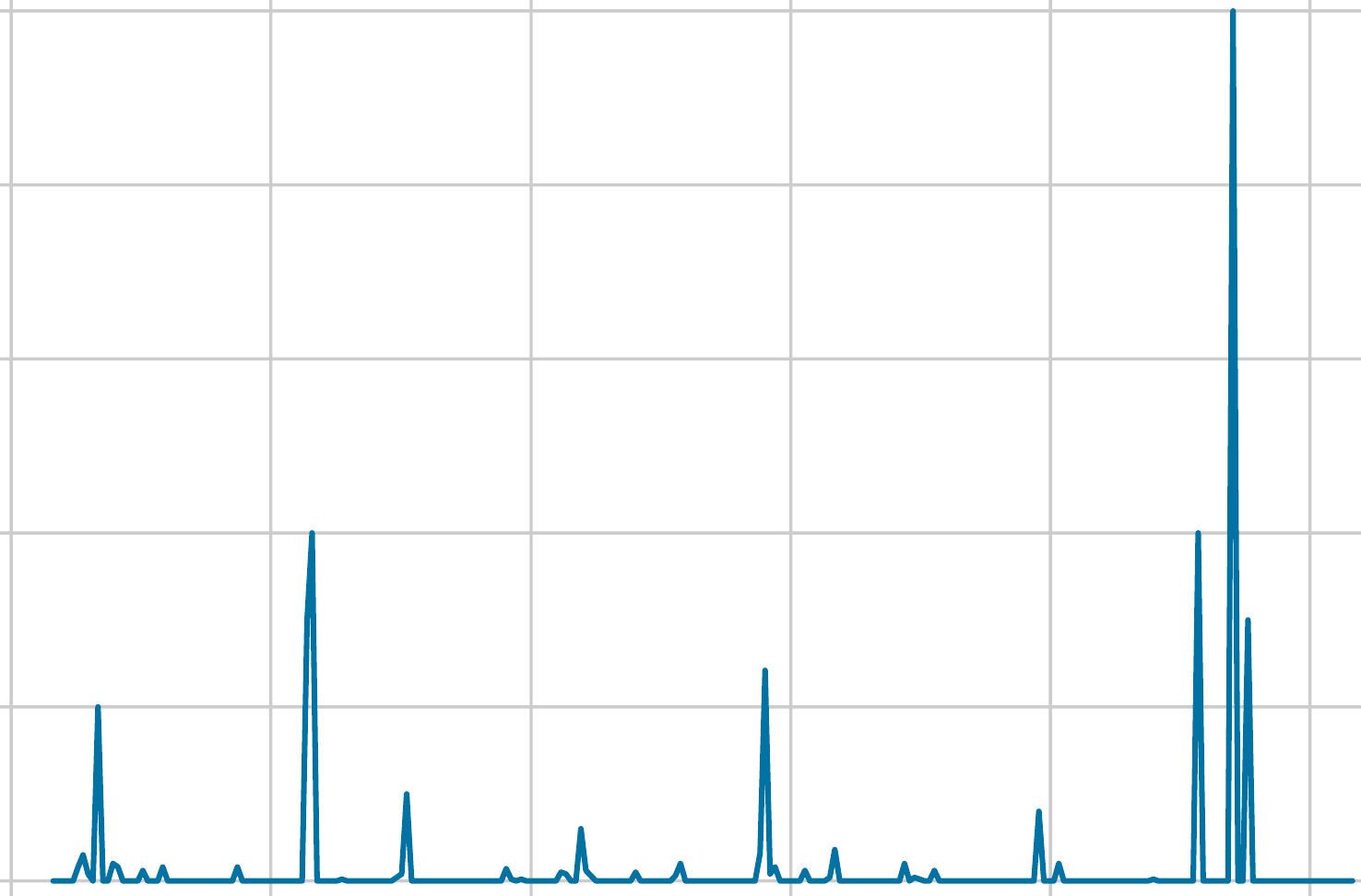}
		\label{fig:ts5}
	}
	\quad
	\subfloat[Monthly aggregation]{
		\includegraphics[width = 3.6cm]{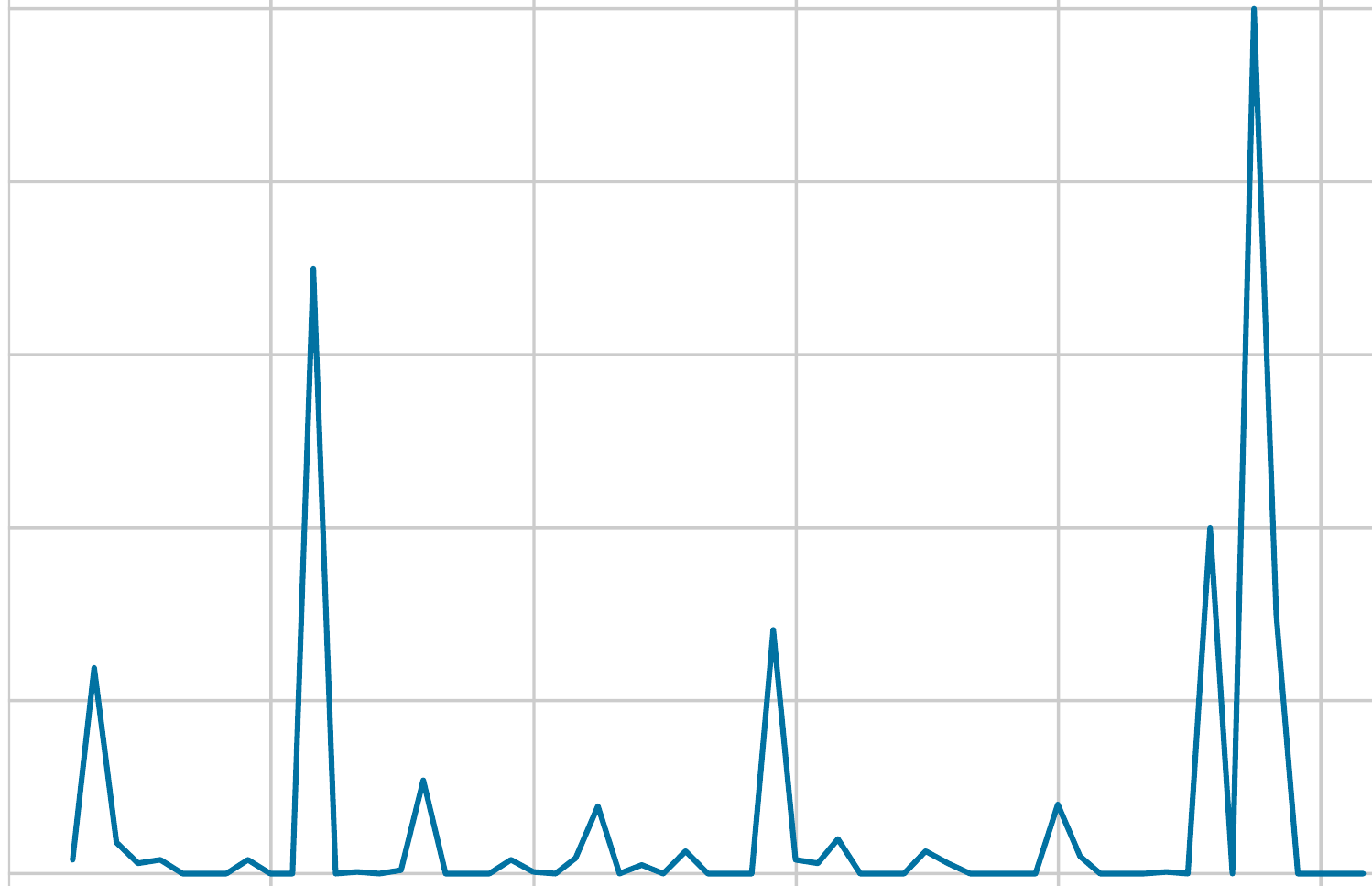}
		\label{fig:ts6}
	}
	\caption{The demand patterns of two SKUs for different demand aggregation levels.}
	\label{fig:bf_time_series}
\end{figure}

The historical demand patterns of SKUs were consolidated based on historical sales and returns. No data imputation was performed, since there were no known missing values in these data sets. The demand patterns of two SKUs are illustrated graphically in Figures~\subref{fig:ts2}--\subref{fig:ts3}, and in Figures~\subref{fig:ts5}--\subref{fig:ts6}, respectively. These two SKUs showcase the diverse characteristics of the historical demand patterns of the SKUs considered (the caption below each figure elucidates the level of aggregation). 

It is clear that the level of intermittency, variation, and skewness in the demand patterns tend to decrease as the aggregation level increases for the first SKU, but remains relatively unchanged for the second. It was therefore deduced that the business may benefit from employing different forecasting models, depending on the historical demand patterns of SKUs. The demand patterns also seemed to contain potential demand exceptions worth considering (something which the business does not currently do).

\subsection{Exception detection}

A number of aggregation levels were considered when detecting potential demand exceptions in order to better understand of the presence of exceptions for the different aggregation levels. The proportion of SKUs for which at least one exception was detected in the historical customer demand patterns, based on the two best-performing DSS methods of the verification analysis, are provided in Table~\ref{tab:Demand_Exceptions_proportion}. 

\begin{table}[pos=h] 
	\caption{Proportion of SKUs with potential demand exceptions.}
	\label{tab:Demand_Exceptions_proportion}
	\begin{tabular}{crr}
		& \multicolumn{2}{c}{Detection method} \\
		\multicolumn{1}{c}{Aggregation level} & \multicolumn{1}{c}{Seasonal} & \multicolumn{1}{c}{AR} \\ \hline
		\multicolumn{1}{c}{Daily} & 100\% & 100\% \\
		\multicolumn{1}{c}{Weekly} & 88.8\% & 99.1\% \\
		\multicolumn{1}{c}{Bi-weekly} & 54.2\% & 89.0\% \\
		\multicolumn{1}{c}{Monthly} & 30.0\% & 62.8\% \\
		\multicolumn{1}{c}{Quarterly} & 4.3\% & 23.4\% \\ \hline
	\end{tabular}%
\end{table}

As expected, the proportion of SKUs exhibiting potential demand exceptions decreased with increased aggregation. A decrease in the proportions of observations identified as demand exceptions per SKU was also observed with increased aggregation. This is important since the business might not want to investigate and select an appropriate handling method (\textit{e.g.}\ Winsorisation) for a large number of observations per SKU, or alter too large a proportion of past demand observations (in the case of automatic exception handling). In the case study, an automatic approach was adopted whereby demand exceptions were reduced to 85\% of their original values. In the remainder of this section, only weekly demand aggregation is considered, since this is the aggregation level adopted by the business.

\subsection{Prioritisation clustering}

Prioritisation clustering was based on four clustering criteria, namely sales volume (measured in units), sales value and profit value (measured in South African Rands), and customer numerousness (\textit{i.e.}\ number of unique customers) for each SKU. These were all assumed to be positively correlated with the economical importance of an SKU and, therefore, important to the long-term financial success of the business. The criteria were extracted for the data forming part of the training set.

Almost all of the variation between SKUs could be attributed to the first three principal components of the data, and clustering was therefore based on these components. A number of $k=2$ clusters was deemed the most appropriate, based on the degree of cluster separation achieved. The $k$-means method performed best in respect of two of the three quality metrics and was employed to generate the final SKU clusters. The mean feature values for the two clusters, given in Table~\ref{tab:BF_sku_seg_values}, were assessed so as to determine the relative importance of the SKUs included in each cluster. 

\begin{table}[pos=h] 
	\caption{Mean feature values of clusters obtained during prioritisation.}
	\label{tab:BF_sku_seg_values}
	\begin{tabular}{crrrc}
		Cluster & \multicolumn{1}{c}{\begin{tabular}[c]{@{}c@{}}Sales\\ volume\end{tabular}} & \multicolumn{1}{c}{\begin{tabular}[c]{@{}c@{}}Sales\\ value\end{tabular}} & \multicolumn{1}{c}{\begin{tabular}[c]{@{}c@{}}Profit\\ value\end{tabular}} & \begin{tabular}[c]{@{}c@{}}Customer\\ numerousness\end{tabular} \\ \hline
		1 $\mapsto$ B & 571\,551 & 364\,672 & 149\,321 & 171 \\
		2 $\mapsto$ A & 5\,207\,002 & 1\,656\,849 & 688\,643 & 741 \\ \hline
	\end{tabular}%
\end{table}

It is clear that the SKUs in Cluster~2 (labelled `A') have historically been more important to the business than SKUs in Cluster~1 (labelled `B'), in terms of all four criteria. Furthermore, Cluster~2 contains only 140 SKUs, while Cluster~1 contains 1125 SKUs. The SKUs in cluster `A,' on average, generated around 910\% higher sales volumes, 450\% higher sales revenue, 460\% more profit, and served 430\% more of the customer base than SKUs in cluster `B'. 

The business does not currently perform SKU classification, other than deciding which SKUs to keep inventory of. Prioritisation clustering in the case study was based on two clusters (which yielded the largest cluster separation) so as to show the potential benefit of prioritising certain SKUs. 

\subsection{Time series clustering and forecasting}

Time series clustering was again based on the features listed in Table~\ref{tab:ts_features}. After performing PCA, the top three principal components were employed. The largest cluster separation was achieved for two clusters. The mean scaled values of the features for which the largest difference between clusters was observed are provided in Table~\ref{tab:ts_cluster_1}. 

\begin{table}[pos=h] 
	\centering
	\caption{Mean scaled values of time series features for clusters.}
	\label{tab:ts_cluster_1}
	\begin{tabular}{lcc}
		& \multicolumn{2}{c}{Cluster} \\
		Feature & 1 & 2 \\ \hline
		Trend & 0.517 & 0.586 \\
		Seasonality & 0.629 & 0.702 \\
		Skewness & 0.447 & 0.192 \\
		Kurtosis & 0.288 & 0.073 \\
		Variation & 0.234 & 0.073 \\
		Intermittency & 0.579 & 0.145 \\ \hline
	\end{tabular}%
\end{table}

The demand patterns of SKUs in Cluster 2 exhibited stronger trend, seasonality, and self-similarity, while the demand patterns of SKUs in Cluster 1 exhibited larger skewness, kurtosis, chaos, variation, and intermittency values. The time series belonging to Cluster~2 are therefore more structured and less noisy, skewed and intermittent. 

Evaluating the difference in mean values (not scaled) for intermittency revealed that, on average, 58.8\% of observations in the demand pattern of Cluster 1 SKUs are equal to zero. For Cluster 2 SKUs, 17.3\% of time series observations are equal to zero. High levels of intermittency were expected, however, based on the observed historical demand patterns.

The time series forecasting models considered during the case study included those of the verification analysis --- the Croston, SBA, and TSB models. During model tuning, forecasts were evaluated based on 5-fold cross validation in respect of the MASE. Evaluating the three best-performing models per time series cluster revealed that the Croston~2, Naive~1, and Naive~2 methods achieved the best cross-validated performance in respect of the Cluster~1 time series, while the Croston~2, SBA, and Theta\_bc models performed the best in respect of the Cluster~2 time series. The difference between performances of the Croston~2 and SBA methods and the other models were statistically significant at $\alpha=0.01$ for the time series in Cluster~2. 

\begin{figure}[pos=h] 
	\centering
	\subfloat[Cluster 1 CD diagram]{
		\includegraphics[scale=0.5]{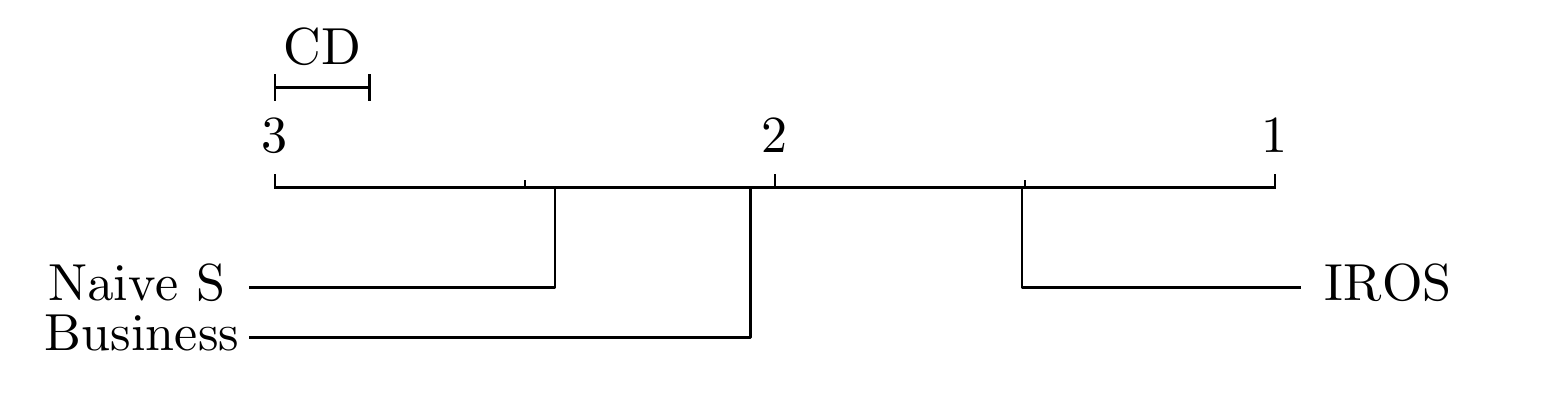}
		\label{fig:BF_weekly_cluster_2_comp_final}
	}
	\quad
	\subfloat[Cluster 2 CD diagram]{
		\includegraphics[scale=0.5]{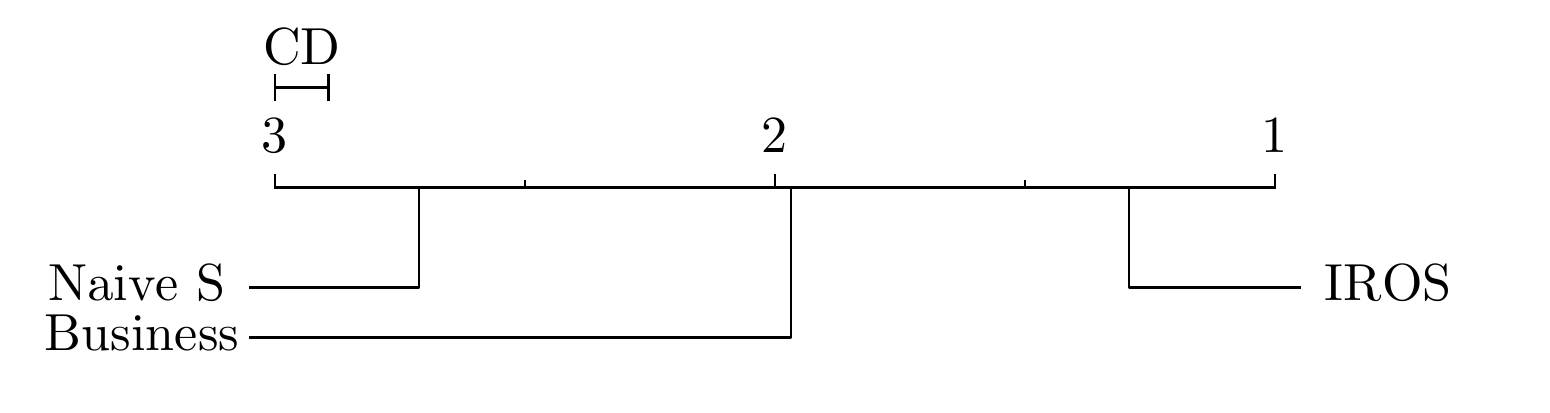}
		\label{fig:BF_weekly_cluster_2_sig_final}
	}
	\caption{Comparison of DSS forecasting with business}
	\label{fig:BF_model_comparison_weekly_final}
\end{figure}

These models were subsequently used during model building (following the same approach as explained in \S\ref{Verification_forecasting}).  When retrospectively comparing the final forecast of the DSS instantiation with the naive mean model currently employed by the business (on the comparison data set), a statistically significant difference in performance was found for $\alpha=0.01$. The CD diagram for each time series cluster, based on the final predictions, is given in Figure~\ref{fig:BF_model_comparison_weekly_final}, where the Naive~S model (which did not perform well during cross validation) was simply included for reference.

For both time series clusters, the DSS instantiation outperformed the current method employed by the business. A larger improvement in forecasting performance was observed for the Cluster~2 time series. This was expected since the historical demand patterns of these time series exhibited stronger patterns which could be captured by more advanced forecasting models. These results also confirmed the suitability of the models specifically developed for intermittent time series when encountering high levels of intermittency. The DSS instantiation outperformed the current forecasting model adopted by the business by 8.8\% (in respect of the MASE) across all the SKUs considered. 

\subsection{Inventory replenishment}

\begin{figure}[pos=h] 
	\centering
	\includegraphics[scale=1.4]{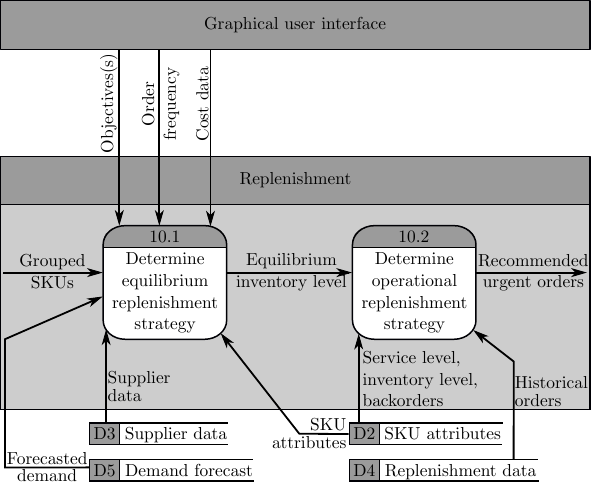}
	\caption{Level-two DFD of Module 10.0 of the framework.}
	\label{fig:replenishment_child}
\end{figure}

The working of the inventory replenishment model implemented in Module~10.0 of the DSS instantiation is depicted by its level-two DFD in Figure~\ref{fig:replenishment_child}. The proposed model is based on two separate optimisation submodels.

An \textit{equilibrium inventory submodel} was first used to determine a replenishment strategy for an ``ideal" equilibrium case. The submodel not only suggests when to place replenishment orders during the planning horizon and how large these orders should be, but also what the starting inventory level and historical orders should be to reach equilibrium (that is, a seasonal repetition of the inventory levels over successive planning horizons). Determining an ideal cyclic inventory level is valuable since it acts as a target towards which the actual (\textit{i.e.}\ operational) inventory level should be managed and is robust to unexpected changes in demand since it does not place restrictions on the manner in which the operational inventory level is managed. 

In a subsequent \textit{operational inventory submodel}, a number of constraints from reality are imposed, such as the actual starting inventory level, backorders, and historical orders. The forecast demand, as well as SKU and container attributes, remain unchanged. Both of these submodels are combinatorial optimisation problems with the same user-specified objective (\textit{e.g.}\ cost minimisation). Both submodels were solved by the branch-and-cut method as implemented in the \texttt{docplex} Python library.

The historical performance of an inventory management system is interpretable both in terms of the overall cost of inventory management, as well as the average inventory levels of SKUs. The business, however, typically places replenishment orders for periods-of-supply which are longer than the anticipated lead time in order to account for late deliveries. The cost of their replenishment approach was determined for varying periods-of-supply and an expected lead time of 12 weeks. The business, however, requested that the specific cost amount remain confidential.

When comparing its cost with the optimised replenishment cost of the DSS instantiation model (with a target service level of 95\% for SKUs in Cluster `A' and 90\% for SKUs in Cluster `B'), over the course of the planning horizon, a potential reduction in replenishment cost of between 11.78\% and 23.53\% was observed, as shown in Table~\ref{fig:Case_study_improvement}. Such savings are significant in the context of the large monetary expenditure and may be attributed to the cost savings of only ordering full container loads in the DSS instantiation model, as well as the fact that the model penalises higher inventory levels.

\begin{table}[pos=h]
	\centering
	\caption{The potential cost savings for the industry partner.}
	\label{fig:Case_study_improvement}
	\begin{tabular}{cc}
		Periods-of-supply & Cost savings (\%) \\ \hline
		24 weeks & 23.53 \\
		20 weeks & 17.49 \\
		16 weeks & 15.56 \\
		12 weeks & 11.78 \\ \hline
	\end{tabular}
\end{table}

In a general retail context, lower inventory levels are usually associated with a decrease in costs related to inventory management. The potential decrease in the mean inventory level of each SKU, when applying the proposed inventory replenishment model, was therefore analysed. The mean decrease in inventory level for all SKUs considered are provided in Table~\ref{fig:Case_study_levels_improvement}. 

\begin{table}[pos=h]
	\centering
	\caption{The potential decrease in inventory levels.}
	\label{fig:Case_study_levels_improvement}
	\begin{tabular}{cc}
		Periods-of-supply & Inventory level decrease (\%) \\ \hline
		24 weeks & 55.63 \\
		20 weeks & 53.38 \\
		16 weeks & 50.49 \\
		12 weeks & 46.86 \\ \hline
	\end{tabular}
\end{table}

Based on these results, the DSS inventory replenishment model realised a mean decrease in inventory level of between 46.86\% and 55.63\%. The benefit of employing a more sophisticated model to determine recommended replenishment orders is even more clear in the inventory level analysis.

When considering the potential cost savings provided in Table~\ref{fig:Case_study_improvement} in conjunction with the potential decrease in inventory levels provided in Table~\ref{fig:Case_study_levels_improvement}, the replenishment component of the IROS framework was judged to have significantly outperformed the current model employed at the case study business.

\section{Discussion}\label{Discussion}

In this paper, a generic framework for decision support in inventory management was proposed. The framework addresses a number of major shortcomings of existing frameworks for inventory management available in the literature. These primarily pertain to the absence of a holistic, integrated approach and the domain-specific nature of many of the frameworks, the vague guidance provided to readers with respect to framework implementation, and the general absence of means toward incorporating the expert knowledge of the user.

\subsection{Research contributions}\label{Contributions}
The IROS framework proposed in this paper is \textit{generic} in nature and its design allows for application in a wide variety of retail business types, whereas many existing frameworks are aimed at domain-specific inventory management problems only. The framework allows for the implementation of simple to state-of-the-art modelling approaches in each of its components, based on user preference.

All of the components one would expect to find in a framework for inventory management are combined in the design of the IROS framework. The IROS framework is also specifically aimed at utilising any expert or domain knowledge of the user by eliciting user input, whenever possible and appropriate. Furthermore, the user is afforded flexibility in the choice of models and methods to employ for the different processes embedded within the framework components. Finally, the IROS framework is presented in an intuitive, easy-to-follow manner, with its processes, and their corresponding inputs and outputs, visualised by means of DFDs.

An approach towards reducing the computational burden of evaluating the cross-validated performance of forecasting models on a large number of time series was also presented. It was shown that the adoption of a time series clustering approach may be beneficial. According to this approach, the performance of different forecasting models are evaluated on a subset of the time series of each cluster and their relative performances are then taken as a proxy of model performance for the entire time series cluster.

\subsection{Future work}\label{Future_work}

The claims pertaining to the general applicability of the IROS framework may further be substantiated by applying it to a case study from an alternative retail domain. The scope of the framework's applicability may be enlarged by, for instance, applying an instantiation of the framework to a case study involving perishable products. The framework may also be adapted so as to support multi-echelon inventory replenishment. 

The implementation of a fully integrated DSS instantiation of the framework by means of a robust, user-friendly graphical user interface, with pre-populated framework modules, may prove to be an effective concept demonstrator when showcasing the capabilities of the IROS framework to inventory management practitioners in industry. Such a computerised proof-of-concept implementation will also aid users in applying the IROS framework. 

The suite of verification data may be expanded to include multivariate time series. With the large volumes of data available to modern retail businesses, exploring the utilisation of these data may well yield favourable results. 

\section{Conclusion}\label{Conclusions}

The resulting output of the IROS framework may prove useful to inventory management practitioners by providing decision support in respect of the main decisions required during inventory replenishment. The application of a DSS instantiation of the framework on benchmark data from the literature showed that the methodology proposed by the design of the framework performed well when compared with other benchmark models. The subsequent application of this verified DSS instantiation, in the context of a case study showed the additional value it may add to the operations of a successful business currently operating in the South African retail sector.

\printcredits

\bibliographystyle{cas-model2-names}

\bibliography{bibliography}

\bio{}

\endbio

\bio{}
\endbio

\end{document}